\documentclass{amsart}
\usepackage{amssymb}

\usepackage{a4}
\begin{document}
 \makeatletter
  \renewcommand{\theequation}{%
   \thesection.\arabic{equation}}
  \@addtoreset{equation}{section}
 \makeatother
\newtheorem{theorem}{Theorem}[section]
\newtheorem{lemma}[theorem]{Lemma}
\newtheorem{corollary}[theorem]{Corollary}
\newtheorem{remark}[theorem]{Remark}
\def\operatorname#1{{\rm#1\,}}
\def\qedbox{\hbox{$\rlap{$\sqcap$}\sqcup$}}
\def\range{\operatorname{range}}
\def\Pspan{\operatorname{span}}
\def\rank{\operatorname{rank}}
\def\id{\operatorname{Id}}
\def\imag{\operatorname{Im}}
\def\Tr{\operatorname{Tr}}
\def\Class{\operatorname{Class}}
\def\MSpin{{\operatorname{MSpin}}}
\def\TKSp{\tilde K Sp}
\def\Tko{\tilde k o}
\def\TMSpin{\tilde M Spin}
\def\TKU{\tilde K U}
\def\KSp{{\operatorname{KSp}}}
\def\ko{{\operatorname{ko}}}
\def\Clif{\operatorname{Clif}}
\def\Irr{\operatorname{Irr}}
\font\pbglie=eufm10
\def\Z{{\mathbb{Z}}}\def\ZL{{\mathbb{Z}_\ell}}
\def\QL{{Q_\ell}}\def\R{{\mathbb{R}}}
\def\Q{{\mathbb{Q}}}\def\HH{{\mathbb H}}\def\C{{\mathbb C}}
\def\fl{{\operatorname{flat}}}
\title[The eta invariant]{The eta invariant and the real connective $K$-theory of the classifying space for quaternion groups}
\author{Egidio Barrera-Yanez and
Peter B. Gilkey}
\begin{address}{Instituto de Matem\'aticas UNAM
U. Cuernavaca, Av. Universidad s/n, Col. Lomas de Chamilpa C.P. 62210,
Cuernavaca, Mor. MEXICO\newline email: ebarrera@matcuer.unam.mx}\end{address}
\begin{address}{Mathematics Dept., University of  Oregon, Eugene, Oregon 97403, USA and\newline
Max Planck Institute for 
Mathematics in the Sciences, 
Inselstrasse 22-26, 04103 Leipzig, Germany. email: gilkey@math.uoregon.edu}\end{address}
\begin{abstract} We express the real connective $K$ theory groups $\Tko_{4k-1}(B\QL)$ of the quaternion
group $\QL$ of order $\ell=2^j\ge8$ in terms of the representation theory of
$\QL$ by showing $\Tko_{4k-1}(B\QL)=\TKSp(S^{4k+3}/\tau\QL)$ where $\tau$ is any fixed point free representation of $\QL$ in $U(2k+2)$.
\end{abstract}
\maketitle

\section{Introduction}\label{Sect-Introduction}

 A compact Riemannian manifold $(M,g)$ is said to be a {\it spherical space form} if $(M,g)$ has constant
sectional curvature $+1$.
A finite group $G$ is said to be a {\it spherical space form group} if there exists a representation
$\tau:G\rightarrow U(k)$ for $k\ge2$ which is fixed point free - i.e. $\det(I-\tau(\xi))\ne0$
$\forall$ $\xi\in G-\{1\}$. Let
$$M^{2k-1}(G,\tau):=S^{2k-1}/\tau(G)$$
be the associated spherical space form; $G$ is then the fundamental group of the manifold $M^{2k-1}(G,\tau)$. Every odd dimensional spherical
space form arises in this manner; the only even dimensional spherical space forms are the sphere $S^{2k}$ and real projective space
$\mathbb{RP}^{2k}$. The spherical space form groups all have periodic cohomology;
conversely, any group with periodic cohomology acts without fixed points on some sphere, although not necessarily orthogonally. We refer to
\cite{Wo85} for further details concerning spherical space form groups.

Any cyclic group is a spherical space form group since the group of $\ell^{th}$ roots of unity acts without fixed points by complex
multiplication on the unit sphere $S^{2k-1}$ in $\C^k$. Let
$\HH=\text{span}_{\mathbb{R}}\{1,\mathcal{I},\mathcal{J},\mathcal{K}\}$ be the quaternions, let $\ell=2^j\ge8$,  and let
$\xi:=e^{4\pi\mathcal{I}/\ell}\in\HH$ be a primitive $(\frac\ell2{})^{th}$ root of unity. The quaternion group
$\QL$ is the subgroup of $\HH$ of order $\ell$ generated by $\xi$ and $\mathcal{J}$:
\begin{equation}
Q_{\ell}:=\{1,\xi,...,\xi^{\ell/2-1},\mathcal{J},\xi\mathcal{J},...,\xi^{\ell/2-1}\mathcal{J}\}.\label{eqn1.1}\end{equation}

Let $BG$ be the classifying space of a finite group and let $ko_*(BG)$ be the associated real connective $K$ theory groups; we refer
to \cite{BaBe00,BaBr96,BG00,GeTh01,Gr99} for a further discussion of connective $K$ theory and related matters. 

The $p$ Sylow subgroup of a spherical space form group $G$ is cyclic if $p$ is odd and either cyclic or a quaternion group $\QL$ for
$\ell=2^j\ge8$ if $p=2$. This focuses attention on these two groups. We showed previously in \cite{BYG} that:
\begin{theorem}\label{theorem1.1}
Let $\ZL$ be the cyclic group of order $\ell=2^j>1$. Let $k\ge1$. Let $\tau:\ZL\rightarrow U(2k+2)$ be a fixed point free
representation. Then
$$\Tko_{4k-1}(B\ZL)=\TKSp(M^{4k+3}(\ZL,\tau)).$$
\end{theorem}

 In this paper, we generalize Theorem \ref{theorem1.1} to the quaternion group:
\begin{theorem}\label{theorem1.2} Let $\QL$ be the quaternion group of order $\ell=2^j\ge3$. Let $k\ge1$. Let
$\tau:\QL\rightarrow U(2k+2)$ be a fixed point free representation. Then 
$$\Tko_{4k-1}(B\QL)=\TKSp(M^{4k+3}(\QL,\tau)).$$
\end{theorem}

The quaternion (symplectic) $K$ theory groups $\TKSp(M^{4k+3}(\QL,\tau))$ are expressible in
terms of the representation theory - see Theorem \ref{theorem4.1}. Thus Theorem \ref{theorem1.2} expresses $\Tko_{4k-1}(B\QL)$ in
terms of representation theory. If $\ell=8$, then these groups were determined previously \cite{BaBr96,BoGi95}.

Here is a brief outline to this paper. In Section
\ref{Sect-RepTheory}, we review some facts concerning the representation theory of
$\QL$ which we shall need. In Section \ref{Sect-Eta}, we review some results concerning the eta invariant. In Section \ref{Sect-KSp},
we use the eta invariant to study $\TKSp(M^{4k+3}(\QL,\tau))$. In Section \ref{Sect-ko}, we use the eta
invariant to study
$\Tko(B\QL)$ and complete the proof of Theorem \ref{theorem1.2}. 

The proof of Theorem \ref{theorem1.2} is quite a bit different from
the proof  of Theorem \ref{theorem1.1} given previously; the extension is not straightforward. This arises from the fact that unlike
the classifying space
$B\ZL$, the $2$ localization of $B\QL$ is not irreducible. Let
$SL_2(\mathbb{F}_q)$ be the group of
$2\times 2$ matrices of determinant
$1$ over the field
$\mathbb{F}_q$ with
$q$ elements where $q$ is odd. Then the $2$-Sylow subgroup of $SL_2(\mathbb{F}_q)$ is $\QL$ for $\ell=2^j$ where $j$ is
the power of $2$ dividing $q^2-1$.  There is a stable 2-local splitting of the classifying space $B\QL$ in the form
\begin{equation}B\QL=BSL_2(\mathbb{F}_q)\vee\Sigma^{-1}BS^3/BN\vee\Sigma^{-1}BS^3/BN\label{eqn4.6}\end{equation}
where $N$ is the normalizer of a maximal torus in $S^3$ \cite{MP84,MP91}. 
It is necessary to find a corresponding splitting of $\TKSp(M^{4k+3}(\QL,\tau))$ that mirrors this
decomposition; see Remark \ref{remark5.2}.

\section{The Representation Theory of $\QL$}\label{Sect-RepTheory}

We say that
$f:\QL\rightarrow\C$ is a {\it class function} if
$f(xgx^{-1})=f(g)$ for all
$x,g\in\QL$; let $\Class(\QL)$ be the Hilbert space of all class functions with the $L^2$ inner product
$$\langle f_1,f_2\rangle=\ell^{-1}\textstyle\sum_{g\in\QL}f_1(g)\bar f_2(g).$$
Let $\Irr(\QL)$ be a set of representatives for the equivalence classes of irreducible unitary representations of $\QL$. The {\it
orthogonality relations} show that
$\{\Tr(\sigma)\}_{\sigma\in\Irr(\QL)}$ is an orthonormal basis for
$\Class(\QL)$, i.e. we may expand any class function:
$$f=\textstyle\sum_{\sigma\in\Irr(\QL)}\langle f,\Tr(\sigma)\rangle\Tr(\sigma).$$

The {\it unitary group representation ring} $RU(\QL)$ and the {\it augmentation ideal} $RU_0(\QL)$ are defined by:
\begin{eqnarray*}
RU(\QL)&=&\text{Span}_\Z\{\sigma\}_{\sigma\in\Irr(\QL)},\text{ and }\\
RU_0(\QL)&=&\{\sigma\in RU(\QL):\dim\sigma=0\}.\end{eqnarray*}
We shall identify a representation with
the class function defined by its trace henceforth; a class function $f$ has the form $f=\Tr(\tau)$ for some $\tau\in RU(\QL)$ if and
only if
$\langle f,\sigma\rangle\in\Z$ for all $\sigma\in\Irr(\QL)$.

Let $RSp(\QL)$ and $RO(\QL)$ be the $\Z$ vector spaces generated by equivalence
classes of irreducible quaternion and real representations, respectively. Forgetting the symplectic structure
and complexification of a real structure define natural inclusions $RSp(\QL)\subset RU(\QL)$ and $RO(\QL)\subset RU(\QL)$. 
We have:
\begin{eqnarray}
&&RO(\QL)\cdot RO(\QL)\subset RO(\QL),\nonumber\\
&&RSp(\QL)\cdot RSp(\QL)\subset RO(\QL),\label{eqn2.1}\\
&&RO(\QL)\cdot RSp(\QL)\subset RSp(\QL).\nonumber
\end{eqnarray}

\medbreak 
The $\frac\ell4+3$ conjugacy classes of $\QL$ have
representatives:
$$\{1,\ \xi,\ ...,\ \xi^{\ell/4}=-1,\ \mathcal{J},\ \mathcal{\xi J}\}.$$
There are $\frac\ell4+3$ irreducible inequivalent complex representations of $\QL$. Four of these representations are the
$1$ dimensional representations defined by:
\begin{eqnarray*}
\begin{array}{llll}
\rho_0(\xi)=1,&\kappa_1(\xi)=-1,& \kappa_2(\xi)=1,&\kappa_3(\xi)=-1,\\
\rho_0(\mathcal{J})=1,&\kappa_1(\mathcal{J})=1,&\kappa_2(\mathcal{J})=-1,&\kappa_3(\mathcal{J})=-1.\end{array}\end{eqnarray*}
We define representations $\gamma_u:\QL\rightarrow  U(2)$ by setting:
\begin{eqnarray*}
\gamma_u(\xi)&=&\left(\begin{array}{rr}\xi^u&0\\0&\xi^{-u}\end{array}\right),\quad
\gamma_u(\mathcal{J})=\left(\begin{array}{rr}0&(-1)^u\\1&0\end{array}\right).\end{eqnarray*}
The representations $\gamma_u$, $\gamma_{-u}$, and
$\gamma_{u+\frac\ell2}$ are all equivalent. The representations $\gamma_u$ are irreducible and inequivalent for 
$1\le u\le\frac\ell4-1$;
$\gamma_0$ is equivalent to $\rho_0+\kappa_2$ and $\gamma_{\frac\ell4}$ is equivalent to $\kappa_1+\kappa_3$. We have:
$$\Irr(\QL)=\{\rho_0,\kappa_1,\kappa_2,\kappa_3,\gamma_1,...,\gamma_{\frac\ell4-1}\}.$$
If $\vec s=(s_1,...,s_k)$ is a $k$ tuple of odd integers, then 
$$\gamma_{\vec s}:=\gamma_{s_1}\oplus...\oplus\gamma_{s_k}$$ 
is a fixed point free
representation from
$\QL$ to $U(2k)$; conversely, every fixed point free representation of $\QL$ is conjugate to such a representation. The associated
spherical space forms are the {\it quaternion spherical space forms}.

The representations
$\{\rho_0,\kappa_1,\kappa_2,\kappa_3\}$ are real, the representations $\gamma_{2i}$ are real, and the representations
$\gamma_{2i+1}$ are quaternion. We have:
$$\begin{array}{l}
RO(\QL)\phantom{.}=\text{span}_\Z\{\phantom{2}\rho_0, 
\phantom{2}\kappa_1,\phantom{2}\kappa_2,\phantom{2}\kappa_3,2\gamma_1,\phantom{2}\gamma_2,...,2\gamma_{\ell/4-1}\},\\
RSp(\QL)=\text{span}_\Z\{2\rho_0,2\kappa_1,2\kappa_2,2\kappa_3,\phantom{2}\gamma_1,2\gamma_2,..., 
\phantom{2}\gamma_{\ell/4-1}\}.\end{array}$$

We define:
\begin{equation}\label{eqn2.2}\begin{array}{ll}
\Theta_1(g)&:=\left\{\begin{array}{ll}
\textstyle\frac\ell4&\text{if }g=\pm\mathcal{I},\\
-2&\text{if }g=\xi^{2i}\mathcal{J},\\
0&\text{ otherwise,}\end{array}\right.\\ \\
\Theta_2(g)&:=\left\{\begin{array}{ll}
\textstyle\frac\ell4&\text{if }g=\pm\mathcal{I},\\
-2&\text{if }g=\xi^{2i+1}\mathcal{J},\\
0&\text{ otherwise.}\end{array}\right.\end{array}\end{equation}
The two class functions $\Theta_i$ will be used to mirror in $RU(\QL)$ the splitting of $B\QL$ given in equation
(\ref{eqn4.6}).

We identify virtual representations with the class functions they define henceforth. Let
$$\Delta:=2\rho_0-\gamma_1;\qquad\Tr(\Delta)=\det(I-\gamma_1).$$

\begin{lemma}\label{lemma2.1}\ \begin{enumerate}
\item We have $\Theta_1\in RO_0(\QL)$ and $\Theta_2\in RO_0(\QL)$.\item Let
$c_i:=\ell^{-1}\textstyle\sum_{g\in\QL-\{1\}}\Delta(g)^{i}$. We have $c_0=\frac{\ell-1}\ell$.
 If $i>0$, then $c_{2i}\in\Z$ and $c_{2i-1}\in2\Z$.
\end{enumerate}\end{lemma}

\medbreak\noindent{\bf Proof:}
We use equation (\ref{eqn2.2}) to compute:
$$\begin{array}{llll}
\text{ for any }\ell&
\langle\Theta_1,\rho_0\rangle=0,&\langle\Theta_1,\gamma_{2i+1}\rangle=0,&\langle\Theta_1,\gamma_{2i}\rangle=(-1)^i,\\
&\langle\Theta_2,\rho_0\rangle=0,&\langle\Theta_2,\gamma_{2i+1}\rangle=0,&\langle\Theta_2,\gamma_{2i}\rangle=(-1)^i,\\
\text{ for }\ell=8&
\langle\Theta_1,\kappa_1\rangle=-1,&\langle\Theta_1,\kappa_2\rangle=1,&\langle\Theta_1,\kappa_3\rangle=0,\\
&\langle\Theta_2,\kappa_1\rangle=0,&\langle\Theta_2,\kappa_2\rangle=1,&\langle\Theta_2,\kappa_3\rangle=-1,\\
\text{ for }\ell>8&
\langle\Theta_1,\kappa_1\rangle=0,&\langle\Theta_1,\kappa_2\rangle=1,&\langle\Theta_1,\kappa_3\rangle=1,\\
&\langle\Theta_2,\kappa_1\rangle=1,&\langle\Theta_2,\kappa_2\rangle=1,&\langle\Theta_2,\kappa_3\rangle=0.\\
\end{array}$$
We use equation (\ref{eqn2.1}) to complete the proof of assertion (1):
$$
\begin{array}{l}
\Theta_1=\left\{\begin{array}{ll}
\Tr\{\kappa_2-\kappa_1\}&\text{if }\ell=8,\\
\Tr\{\kappa_2+\kappa_3+\textstyle\sum_{1\le i<\ell/8}(-1)^i\gamma_{2i}\}&\text{if }\ell\ge16,
\end{array}\right.\\
\Theta_2=\left\{\begin{array}{ll}
\Tr\{\kappa_2-\kappa_3\}&\text{if }\ell=8,\\
\Tr\{\kappa_2+\kappa_1+\textstyle\sum_{1\le i<\ell/8}(-1)^i\gamma_{2i}\}&\text{if }\ell\ge16.
\end{array}\right.\end{array}
$$

The first identity of assertion (2) is immediate. Let $r>0$. As $\Tr(\Delta^r)(1)=0$,
$$c_r=\ell^{-1}\textstyle\sum_{g\in\QL-\{1\}}\Tr(\Delta^r)(g)=\langle\Delta^r,\rho_0\rangle\in\Z.$$
If $r$ is odd, then $\gamma_1^r$ is quaternion so $\langle\gamma_1^r,\rho_0\rangle\in2\Z$. Since
$\Delta^r\equiv\gamma_1^r$ mod $2RU(\QL)$, $\langle\Delta^r,\rho_0\rangle\in2\Z$ if $r$ is odd.
\qedbox

\section{The eta invariant, $K$ theory, and bordism}\label{Sect-Eta}

Let $V$ be a smooth complex vector bundle over a compact Riemannian manifold $M$. Let $V$ be equipped with a unitary
(Hermitian) inner product. Let 
$$P:C^\infty(V)\rightarrow C^\infty(V)$$
be a self-adjoint elliptic first order partial differential operator. Let
$\{\lambda_i\}$ denote the eigenvalues of $P$ repeated according to multiplicity. Let
$$\eta(s,P):=\textstyle\sum_i\operatorname{sign}(\lambda_i)|\lambda_i|^{-s}.$$
The series defining $\eta$ converges absolutely for $\Re(s)>>0$ to define a holomorphic function of $s$. This
function has a meromorphic extension to the entire complex plane with isolated simple poles. The value
$s=0$ is regular and one defines
$$\eta(P):=\textstyle\frac12\{\eta(s,P)+\dim(\ker P)\}|_{s=0}$$
as a measure of the spectral asymmetry of $P$; we refer to \cite{PG95} for further details concerning this invariant which was first
introduced by \cite{APS} and which plays an important role in the index theorem for manifolds with
boundary.

We say that $P$ is {\it quaternion} if $V$ has a quaternion structure and if the action of $P$ commutes with
this structure. We say that $P$ is {\it real} if $V$ is the complexification of an underlying real vector bundle and if $P$ is
the complexification of an underlying real operator.\goodbreak

\begin{lemma}\label{lemma3.1} Let $M$ be a spin manifold of dimension $m$.
\begin{enumerate}
\item If $m\equiv3,4$ mod $8$, then the Dirac operator is quaternion.
\item If $m\equiv7,8$ mod $8$, then the Dirac operator is real.
\end{enumerate}\end{lemma}

\medbreak\noindent{\bf Proof:} Let $\Clif(m)$ be the real Clifford algebra on $\R^m$. We have:
$$\begin{array}{l}
\Clif(3)=\HH\oplus\HH,\\
\Clif(4)=M_2(\HH),\\
\Clif(7)=M_8(\R)\oplus M_8(\R),\\
\Clif(8)=M_{16}(\R),\text{ and}\\
\Clif(m+8)=\Clif(m)\otimes_\R M_{16}(\R).\end{array}$$
Therefore, the fundamental spinor representation of $\Clif(m)$ is quaternion if we have $m\equiv3,4$ mod $8$ and real if we have
$m\equiv7,8$ mod $8$. The Lemma now follows. \qedbox 

\medbreak The following deformation result will be crucial to our investigations:

\begin{lemma}\label{lemma3.2} Let $P_u$ be a smooth $1$ parameter family of self-adjoint
first order elliptic partial differential operators on a compact manifold $M$. \begin{enumerate}
\item The reduction mod $\Z$ of
$\eta(P_u)$ is a smooth $\R/\Z$ valued function.
\item The variation $\partial_u\eta(P_u)$ is locally computable. 
\item If the operators
$P_u$ are quaternion, then the reduction mod $2\Z$ of $\eta(P_u)$ is a smooth $R/2\Z$ valued function.
\end{enumerate}
\end{lemma}

\medbreak\noindent{\bf Proof:} We sketch the proof briefly and refer to \cite{PG95} Theorem 1.13.2 for further details. Since
$\frac12\operatorname{sign}(u)$ has an integer jump when $u=0$, 
$\eta(P_u)$ can have integer valued jumps at values of $u$ where $\dim(\ker(P_u))>0$. However, in $\R/\Z$, the jump disapears so the
mod
$\Z$ reduction of
$\eta(P_u)$ is a smooth $\R/\Z$ valued function of $u$; one uses the pseudo-differential calculus to construct an approximate
resolvant and to show that the variation
$\partial_u\eta(P_u)$ is locally computable. Assertions (1) and (2) then follow. If $P_u$ is quaternion, then the eigenspaces
of
$P_u$ inherit quaternion structures. Thus $\dim(\ker P_u)$ is even so $\eta(P_u)$ has twice integer jumps as
eigenvalues cross the origin. Consequently the reduction mod $2\Z$ of $\eta(P_u)$ is smooth and
assertion (3) follows.
\qedbox

\medbreak  

Let $\tilde M$ be the universal cover of a connected manifold $M$ and let $\sigma$ be a representation of $\pi_1(M)$ in $U(k)$. The
associated vector bundle is defined by:
\begin{eqnarray*}&&V^\sigma:=\tilde M\times\C^k/\sim\text{ where we identify}\\
&&(\tilde x,z)\sim(g\cdot\tilde x,\sigma(g)\cdot z)\text{ for }g\in\pi_1(M),\ \tilde x\in\tilde M,\text{ and }z\in\C^k.
\end{eqnarray*}
The trivial connection on $\tilde M\times\C^k$ descends to define a flat connection on $V^\sigma$. The transition functions
of $V^\sigma$ are locally constant; they are given by the representation $\sigma$. Thus the bundle $V^\sigma$ is said to be {\it
locally flat}. Let $P:C^\infty(V)\rightarrow C^\infty(V)$ be a self-adjoint elliptic first order operator on $M$;
$$P^\sigma:C^\infty(V\otimes V^\sigma)\rightarrow C^\infty(V\otimes V^\sigma)$$
is a well defined operator which is locally isomorphic to $k$ copies of $P$.  Define $\eta^\sigma(P):=\eta(P^\sigma)$; we extend by
linearity to $\sigma\in RU(\pi_1(M))$. 

\medbreak This invariant is a homotopy invariant.

\begin{lemma}\label{lemma3.3} Let $P_u$ be a smooth $1$ parameter family of elliptic first order self-adjoint partial differential
operators over $M$.
\begin{enumerate}
\item If $\sigma\in RU_0(\pi_1(M))$, then the mod $\Z$ reduction of $\eta^\sigma(P_u)$ is
independent of the parameter $u$.  
\item If all the operators $P_u$ are quaternion and $\sigma\in RO_0(\pi_1(M))$ or  if all the operators $P_u$
are real and
$\sigma\in RSp_0(\pi_1(M))$, then the mod $2\Z$ reduction of $\eta(P_u,\sigma)$ is independent of the parameter $u$.
\end{enumerate}\end{lemma}

\medbreak\noindent{\bf Proof:} If $\sigma$ is a representation of $\pi_1(M)$, then the mod $\Z$ reduction of $\eta^\sigma(P_u)$ is
smooth a smooth function of $u$ by Lemma \ref{lemma3.2}. Since
$P^\sigma_u$ is locally isomorphic to $\dim\sigma$ copies of $P_u$ and since the variation is locally computable,
$$\partial_u\eta^\sigma(P_u)=\dim\sigma\cdot\partial_u\eta(P_u).$$
This formula continues to hold for virtual representations. In particular, if we have that $\sigma\in RU_0(\pi_1(M))$,
then $\dim\sigma=0$ so $\partial_u\eta^\sigma(P_u)=0$; (1) follows. 

If $P_u$ is quaternion and $\sigma$ is real or if
$P_u$ is real and if
$\sigma$ is quaternion, then $P_u^\sigma$ is quaternion and $\eta^\sigma(P_u)$ is a smooth $\R/2\Z$ valued function of $u$.
The same argument shows that $\partial_u\eta^\sigma(P_u)=0$. \qedbox

\medbreak We can use the eta invariant to construct invariants of $K$ theory.  
Let $P:C^\infty(V)\rightarrow C^\infty(V)$ be a first order self-adjoint elliptic partial
differential operator with leading symbol $p$. Let
$W$ be a unitary vector bundle over
$M$. We use a partition of unity to construct a self-adjoint elliptic first order operator
$P^W$ on $C^\infty(V\otimes W)$ with leading symbol $p\otimes\operatorname{id}$; this operator is not, of course, cannonically
defined.

We can extend the invariant $\eta^\sigma$ to the the reduced unitary unitary and quaternion (symplectic) $K$
theory groups $\TKU$ and $\TKSp$:

\begin{theorem}\label{theorem3.4} Let $P$ be an elliptic self-adjoint first order partial differential
operator. Let $\sigma\in RU_0(\pi_1(M))$.\begin{enumerate}
\item The map $W\rightarrow\eta^\sigma(P^W)$ extends to a map $\eta_P^\sigma:\TKU(M)\rightarrow\R/\Z$.
\item Suppose that $P$ and $\sigma$ are both real or that $P$ and $\sigma$ are both quaternion. 
The map $W\rightarrow\eta^\sigma(P^W)$ extends to a map $$\eta_P^\sigma:\TKSp(M)\rightarrow\R/2\Z.$$
\end{enumerate}
\end{theorem}

\medbreak\noindent{\bf Proof:} Let $P^W$ and $\tilde P^W$ be two first order self-adjoint partial differential operators on
$C^\infty(V\otimes W)$ with leading symbol
$p\otimes\operatorname{id}$. Set:
$$P_u:=uP^W+(1-u)\tilde P^W.$$ 
This is a smooth $1$ parameter family of first order self-adjoint partial differential operators. As the leading symbol of $P_u$
is $p\otimes\operatorname{id}$, the operators $P_u$ are elliptic. By Lemma
\ref{lemma3.3},
$\eta^\sigma(P_u)\in\R/\Z$ is independent of $u$. Consequently
$\eta^\sigma_P(W):=\eta^\sigma(P^W)\in\R/\Z$ only depends on the isomorphism class of the
bundle
$W$. As the eta invariant is additive with respect to direct sums, we may extend $\eta^\sigma_P$ to
$\TKU(M)$ as an $\R/\Z$ valued invariant. Let
$W$ be quaternion. By Lemma \ref{lemma3.3}, $\eta^\sigma(P_u)\in\R/2\Z$ is independent of $u$
if both
$P$ and $\sigma$ are real or if both $P$ and $\sigma$ are quaternion and thus $\eta^\sigma$ extends to $\TKSp$ as an $\R/2\Z$
valued invariant in this instance.
\qedbox

\medbreak We can use the Atyiah-Patodi-Singer index theorem \cite{APS} to see that the eta invariant also defines bordism invariants.
Let
$G$ be a finite group. A $G$ structure $f$ on a connected manifold $M$ is a representation $f$ from $\pi_1(M)$ to $G$. Equivalently,
$f$ can also be regarded as a map from $M$ to the classifying space $BG$. We consider tuples $(M,g,s,f)$
where
$(M,g)$ is a compact Riemannian manifold with a spin structure $s$ and a $G$ structure $f$. We
introduce the bordism relation $[(M,g,s,f)]=0$ if there exists a compact manifold $N$ with boundary $M$ so that the structures
$(g,s,f)$ extend over $N$; this induces an equivalence relation and the equivariant bordism groups $\MSpin_m(BG)$ consists of bordism
classes of these triples. Disjoint union defines the group structure. 

Let $\MSpin_*:=\MSpin_*(B\{1\})$ be defined by the trivial group. Cartesian
product makes $\MSpin_*(BG)$ into an $\MSpin_*$ module.  Let $\mathcal{F}$ be 
the forgetful homomorphism which forgets the $G$
structure $f$. The reduced bordism groups are then defined by:
$$\TMSpin_*(BG):=\ker(\mathcal{F}):\MSpin_*(BG)\rightarrow\MSpin_*.$$
Since the eta invariant vanishes on $\MSpin_*$, we restrict henceforth to the reduced groups.

If $s$ is a spin structure on $(M,g)$, let
$P_{(M,g,s)}$ be the associated Dirac operator. If $\sigma\in
RU_0(G)$, then
$f^*\sigma\in RU_0(\pi_1(M))$ and we may define:
$$\eta^\sigma(M,g,s,f):=\eta^{f^*\sigma}(P_{(M,g,s)}).$$

\begin{theorem}\label{theorem3.5} Let $G$ be a finite group. Assume either that $m\equiv3$ mod $8$ and that
$\sigma\in RO_0(G)$ or that $m\equiv7$ mod $8$ and that $\sigma\in RSp_0(G)$. Then the map $(M,g,s,f)\rightarrow\eta^\sigma(M,s,f)$
extends to a map 
$$\eta^\sigma:\TMSpin_m(BG)\rightarrow\R/2\Z.$$
\end{theorem}

\medbreak\noindent{\bf Proof:} We sketch the proof and refer to \cite{BGS97} for further details. Suppose that
$m\equiv3$ mod $4$ and that $[(M,g,s,f)]=0$ in
$\MSpin_m(BG)$. Then $M=dN$ where the spin and $G$ structures on $M$ extend over $N$. We may also extend the given Riemannian metric
on $M$ to a Riemannian metric on $N$ which is product near the boundary. 

Let
$\sigma\in RU_0(G)$. The Dirac operator
$P_{(M,g,s)}$ on $M$ is the tangential operator of the spin complex
$Q_{(N,g,s)}$ on
$N$. We twist these operators by taking coefficients in the locally flat virtual bundle $V^{f^*\sigma}$.

Let $\hat A(N,g,s)$ be the $A$-roof genus and let
$ch(V^{f^*\sigma})$ be the Chern character. By the Atiyah-Patodi-Singer index theorem \cite{APS}:
$$\operatorname{index}(Q^{f^*\sigma}_{(N,g,s)})=\textstyle\int_N\hat A(N,g,s)\wedge ch(V^{f^*\sigma})+\eta(P^\sigma_{(M,g,s)}).$$
Since $V^{f^*\sigma}$ is a virtual bundle of virtual dimension $0$ which admits a flat connection,
the Chern character of $V^\sigma$ vanishes. Consequently:
$$\eta^\sigma(M,g,s,f)=\eta(P^{f^*\sigma}_{(M,g,s)})=\operatorname{index}(Q^{f^*\sigma}_{N,g,s}).$$

The dimension of $N$ is $m+1$. We apply Lemma \ref{lemma3.1} to see that if $m\equiv3$ mod $8$ and if $\sigma$ is real or if
$m\equiv7$ mod $8$ and if
$\sigma$ is quaternion, then $Q^{f^*\sigma}_{(N,s,f)}$ is quaternion. Thus $\operatorname{index}(Q^{f^*\sigma}_{(N,s,f)})\in2\Z$ so
$\eta^\sigma(M,g,s)$ vanishes as an $\R/2\Z$ valued invariant if $[(M,g,s,f)]=0$ in $\MSpin_m(BG)$.
\qedbox

\medbreak There is a geometric description of the real connective $K$ theory groups $\Tko_m(BG)$ in terms of the spin bordism
groups. Let $\mathbb{HP}^2$ be the quaternionic projective plane. Let $\tilde T_m(BG)$ be the subgroup of $\TMSpin_m(BG)$ consisting of
bordism classes $[(E,g,s,f)]$ where $E$ is the total space of a geometrical $\mathbb{HP}^2$ spin fibration and where the $G$
structure on
$E$ is induced from a corresponding $G$ structure on the base. The following theorem is a special case of a more general result \cite{St94}:

\begin{theorem}\label{theorem3.6} Let $G$ be a finite group. There is a $2$ local isomorphism between $\Tko_m(BG)$ and
$\TMSpin_m(BG)/\tilde T_m(BG)$.
\end{theorem}

We use Theorem \ref{theorem3.6} to draw the following consequence:

\begin{corollary}\label{corollary3.7} Assume either that $m\equiv3$ mod $8$ and
$\sigma\in RO_0(\QL)$ or that $m\equiv7$ mod $8$ and  $\sigma\in RSp_0(\QL)$. Then $\eta^\sigma$ extends to a map
from $\Tko_m(B\QL)$ to $\Q/2\Z$.
\end{corollary}

\medbreak\noindent{\bf Proof:} If $[(E,s,f)]\in T_m(B\QL)$, then
$\eta^\sigma(P_{(E,g,s)})=0$; see \cite{BGS97} Lemma 4.3 or
\cite{refGLP} Lemma 2.7.10 for details. Thus by Theorems \ref{theorem3.5} and Theorem \ref{theorem3.6}, the eta invariant extends
to $\Tko(B\QL)$. By
\cite{BGS97} Theorem 2.4,
$\Tko_{4k-1}(B\QL)$ is a finite $2$ group. Thus it is not necessary to localize at the prime $2$ and the eta invariant takes values
in $\Q/2\Z$. \qedbox

\medbreak The eta invariant is combinatorially computable for spherical space forms. The following theorem follows from \cite{Do78}.

\begin{theorem}\label{theorem3.8} Let
$\tau:G\rightarrow SU(2k)$ be fixed point free, let $P$ be the Dirac operator on $M^{4k-1}(G,\tau)$, and let $\sigma\in
RU_0(G)$. Then
$$\eta^\sigma(P)=\ell^{-1}\textstyle\sum_{g\in G-\{1\}}\Tr(\sigma(g))\det(I-\tau(g))^{-1}.$$
\end{theorem}

\section{The groups $\TKSp(M^{4\nu-1}(\QL,\nu\cdot\gamma_1))$}\label{Sect-KSp}

Let
$\Delta=\det(I-\gamma_1)\in RSp_0(\QL)$. By equation (\ref{eqn2.1}):
$$\begin{array}{ll}
\Delta^\nu RSp(\QL)\subset RSp_0(\QL)&\text{if }\nu\text{ is even},\\
\Delta^\nu RO(\QL)\subset RSp_0(\QL)&\text{if }\nu\text{ is odd}.
\end{array}$$
The following
Theorem is well known - see, for example
\cite{PG89,GK87}:

\begin{theorem}\label{theorem4.1} Let $\tau:\QL\rightarrow U(2\nu)$ be fixed point free. Then
$$\TKSp(M^{4\nu-1}(\QL,\tau))=\left\{
\begin{array}{ll}
RSp_0(\QL)/\Delta^\nu RSp(\QL)&\text{if }\nu\text{ is even},\\
RSp_0(\QL)/\Delta^\nu RO(\QL)&\text{if }\nu\text{ is odd}.\end{array}\right.$$\end{theorem}

By Theorem \ref{theorem4.1}, the particular representation $\tau$ plays no role and we therefore set $\tau=\nu\cdot\gamma_1$.
We use the eta invariant to study these groups. Let $\eta_\nu^\sigma(W)$ be
the invariant described in Theorem \ref{theorem3.4} for the Dirac operator $P$ on $M^{4\nu-1}(\QL,\nu\cdot\gamma_1)$. We define:
$$
\vec\eta_\nu(W):=\left\{\begin{array}{ll}
(\eta_\nu^{\Theta_1\phantom{2}},\eta_\nu^{\Theta_2\phantom{2}},
     \eta_\nu^{2\Delta},\eta_\nu^{\Delta^2\phantom{2}},...,\eta_\nu^{\Delta^{\nu -2}},
\eta_\nu^{2\Delta^{\nu -1}})(W)&\text{if }\nu\text{ is even},\\
(\eta_\nu^{2\Theta_1},\eta_\nu^{2\Theta_2},
     \eta_\nu^{\Delta\phantom{2}},\eta_\nu^{2\Delta^2},...,
\eta_\nu^{\Delta^{\nu -2}},\eta_\nu^{2\Delta^{\nu -1}})(W)&\text{if }\nu\text{ is odd}.
\end{array}\right.$$

\begin{lemma}\label{lemma4.2} Let $M:=M^{4\nu-1}(\QL,\nu\cdot\gamma_1)$. Then $$\vec\eta_\nu:\TKSp(M)\rightarrow(\Q/2\Z)^{\nu +1}.$$
\end{lemma}

\medbreak\noindent{\bf Proof:} We apply Lemma
\ref{lemma3.1} and Theorem \ref{theorem3.4}. We distinguish two cases:
\begin{enumerate}
\item If $\nu$ is even, then $P$ is real. Thus $\eta_\nu^\sigma:\TKSp(M)\rightarrow\Q/2\Z$ for real $\sigma$ and the
Lemma follows as we have used the real representations\newline $\{\Theta_1,\Theta_2,2\Delta,\Delta^2,...,\Delta^{\nu -2},2\Delta^{\nu -1}\}$
to define
$\vec\eta_\nu$.
\item If $\nu$ is is odd, then $P$ is quaternion. Thus
$\eta_\nu^\sigma:\TKSp(M)\rightarrow\Q/2\Z$ if $\sigma$ is quaternion and the Lemma follows as we have used the quaternion
representations
$\{2\Theta_1,2\Theta_2,\Delta,2\Delta^2,...,\Delta^{\nu -2},2\Delta^{\nu -1}\}$ to define $\vec\eta_\nu$. \qedbox\end{enumerate}

\medbreak Let $\varepsilon_{2i}=2$ and $\varepsilon_{2i-1}=1$;
$\{2\Theta_1,2\Theta_2,\Delta,2\Delta^2,...,\varepsilon_{\nu-1}\Delta^{\nu -1}\}$ are quaternion. In Lemma \ref{lemma2.1}, we defined constants
$$c_i:=\ell^{-1}\textstyle\sum_{g\in\QL-\{1\}}\det(I-\gamma_1(g))^i.$$
Since $\Delta(g)=\det(I-\gamma_1(g))$, we use Theorem \ref{theorem3.8} to compute:
\begin{equation}\eta_\nu^{\Delta^r}(\Delta^s)=\ell^{-1}\textstyle\sum_{g\in\QL-\{1\}}\Delta(g)^{r+s}\Delta(g)^{-\nu }=c_{r+s-\nu }.
\label{eqn4x}\end{equation}
Since $\Theta_1$ and $\Theta_2$ are supported on the elements of order $4$ in $Q_\ell$ and since $\Delta(g)=2$ for such an element, we
may use Theorem \ref{theorem3.8} and equation (\ref{eqn2.2}) to see:
\begin{equation}\begin{array}{l}
\eta_\nu^{\Delta^r}(\Theta_i)\phantom{.}=\eta_\nu^{\Theta_i}(\Delta^r)=\ell^{-1}\textstyle\sum_{g\in\QL-\{1\}}
2^r\Tr(\Theta_i(g))2^{-\nu}\\
\qquad\phantom{......}\vphantom{\vrule height12pt}=\ell^{-1}2^{r-\nu}\textstyle\sum_{g\in\QL-\{1\}}\Tr(\Theta_i(g))=0,\\
\vphantom{\vrule height12pt}\eta_\nu^{\Theta_1}(\Theta_1)=\eta_\nu^{\Theta_2}(\Theta_2)=\ell^{-1}2^{-\nu}\textstyle\sum_{g\in\QL-\{1\}}
\Tr(\Theta_1(g))^2\\
\qquad\phantom{......}\vphantom{\vrule height12pt}=\ell^{-1}2^{-\nu}\{2\cdot\textstyle\frac{\ell^2}{16}+4\cdot\frac\ell4\},\\
\vphantom{\vrule height12pt}
\eta_\nu^{\Theta_1}(\Theta_2)=\eta_\nu^{\Theta_2}(\Theta_1)=\ell^{-1}2^{-\nu}\textstyle\sum_{g\in\QL-\{1\}}
\Tr(\Theta_1(g))\Tr(\Theta_2(g))\\
\qquad\phantom{......}\vphantom{\vrule height12pt}=\ell^{-1}2^{-\nu}\{2\cdot\textstyle\frac{\ell^2}{16}\}.
\end{array}\label{eqn4z}\end{equation}
We have $\ell=2^j$. We use equation (\ref{eqn4x}), equation (\ref{eqn4z}), and Lemma \ref{lemma2.1} to see:
$$\vec\eta_\nu\left(\begin{array}{c}
2\Theta_1\\2\Theta_2\\\Delta\\2\Delta^2\\...\\\varepsilon_{\nu-1}\Delta^{\nu -1}\end{array}\right)
=\left(\begin{array}{ll}
A_\nu&0\\0&B_\nu
\end{array}\right)\in M_{\nu+1}(\Q/2\Z)$$
where $A$ is the $2\times2$ matrix given by
$$A_\nu=2^{1-\nu}\left(\begin{array}{cc}
2^{j-3}+1&2^{j-3}\\2^{j-3}&2^{j-3}+1
\end{array}\right)\text{ if }\nu\text{ is even}$$
$$A_\nu=2^{2-\nu}\left(\begin{array}{cc}
2^{j-3}+1&2^{j-3}\\2^{j-3}&2^{j-3}+1
\end{array}\right)\text{ if }\nu\text{ is odd}$$
and where $B$ is the $\nu-1\times \nu-1$ matrix given by:
$$B_\nu=\left(\begin{array}{rrrrrrrr}
2c_{2-\nu }&c_{3-\nu }&2c_{4-\nu }&...&2c_{-2}&c_{-1}&2c_0\\
4c_{3-\nu }&2c_{4-\nu }&4c_{5-\nu }&...&4c_{-1}&2c_0&0\\
2c_{4-\nu }&c_{5-\nu }&2c_{6-\nu }&...&2c_0&0&0\\
...&...&...&...&...&...&...\\
2c_{-2}&c_{-1}&2c_0&...&0&0&0\\
4c_{-1}&2c_0&0&...&0&0&0\\
2c_0&0&0&...&0&0&0
\end{array}\right)\text{ if }\nu\text{ is even}$$
$$B_\nu=\left(\begin{array}{llrrrrrrrr}
c_{2-\nu }&2c_{3-\nu }&c_{4-\nu }&...&2c_{-2}&c_{-1}&2c_0\\
2c_{3-\nu }&4c_{4-\nu }&2c_{5-\nu }&...&4c_{-1}&2c_0&0\\
c_{4-\nu }&2c_{5-\nu }&c_{6-\nu }&...&2c_0&0&0\\
...&...&...&...&...&...&...\\
2c_{-2}&4c_{-1}&2c_0&...&0&0&0\\
c_{-1}&2c_0&0&...&0&0&0\\
2c_0&0&0&...&0&0&0
\end{array}\right)\text{ if }\nu\text{ is odd}.$$

\begin{theorem}\label{theorem4.3} Let $\mathcal{B}_\nu$ be the subgroup of $(\Q/2\Z)^{\nu -1}$ spanned by the rows of the matrix $B_\nu$
defined above. Let $M=M^{4\nu-1}(\QL,\nu\cdot\gamma_1)$. Then 
$$\TKSp(M)=\left\{\begin{array}{ll}
\Z_{2^{\nu }}\oplus\Z_{2^{\nu }}\oplus\mathcal{B}_\nu&\text{if }\nu\text{ is even,}\\
\Z_{2^{\nu -1}}\oplus\Z_{2^{\nu -1}}\oplus\mathcal{B}_\nu&\text{if }\nu\text{ is odd.}
\end{array}\right.$$
\end{theorem}

\medbreak\noindent{\bf Proof:} Let $\mathcal{K}_\nu$ be the subspace of $\TKSp(M)$ spanned by the virtual vector bundles
defined by $\{2\Theta_1,2\Theta_2,\Delta,2\Delta^2,...,\varepsilon_{\nu-1}\Delta^{\nu -1}\}$. It is then immediate from the definition
and from the form of the matrix $A_\nu$ that
\begin{equation}\label{eqn4.1}\vec\eta_\nu(\mathcal{K}_\nu)=\left\{\begin{array}{ll}
\Z_{2^{\nu }}\oplus\Z_{2^{\nu }}\oplus\mathcal{B}_\nu&\text{if }\nu\text{ is even,}\\
\Z_{2^{\nu -1}}\oplus\Z_{2^{\nu -1}}\oplus\mathcal{B}_\nu&\text{if }\nu\text{ is odd.}
\end{array}\right.\end{equation}
We use Lemma \ref{lemma2.1} to see $c_0=\textstyle\frac{\ell-1}\ell$. Thus $2c_0$ is an element of order $\ell$ in $\Q/2\Z$. We use
the diagonal nature of matrix $B_\nu$ to see that:
\begin{equation}\label{eqn4.2}
|\vec\eta_\nu(\mathcal{K}_\nu)|\ge\left\{\begin{array}{ll}
4^{\nu }\ell^{\nu -1}&\text{if }\nu\text{ is even},\\
4^{\nu -1}\ell^{\nu -1}&\text{if }\nu\text{ is odd}.
\end{array}\right.
\end{equation}

The $E_2$ term in the Atiyah-Hirezbruch spectral sequence for the $K$ theory groups $\TKSp^*(M)$ is
$$\oplus_{u+v=w}\tilde{H}^u(M;\KSp^v(pt)).$$
We take $w=0$ and study the reduced groups to obtain the estimate:
\begin{equation}\label{eqn4.3}
|\TKSp(M)|\le|\oplus_{u+v=0}{\tilde H}^u(M;\KSp^v(pt))|.\end{equation}
We have that:
\begin{equation}\label{eqn4.4}\begin{array}{llrl}
\KSp^v(pt)&=&\Z&\text{if }v\equiv0,4\text{ mod }8,\\
\KSp^v(pt)&=&\Z_2&\text{if }v\equiv-5,-6\text{ mod }8,\\
\KSp^v(pt)&=&0&\text{otherwise},\\
\tilde H^u(M;\Z)&=&\ZL&\text{if }u\equiv0,4\text{ mod }8,\ u<4\nu-1,\\
\tilde H^u(M;\Z_2)&=&
\Z_2\oplus\Z_2&\text{if }u\equiv1,2,5,6\text{ mod }8,\ u\le 4\nu-1.
\end{array}\end{equation}
Equations (\ref{eqn4.3}) and (\ref{eqn4.4}) then imply:
\begin{equation}\label{eqn4.5}
|\TKSp(M)|\le\left\{\begin{array}{ll}
4^{\nu }\ell^{\nu -1}&\text{if }\nu\text{ is even}\\
4^{\nu -1}\ell^{\nu -1}&\text{if }\nu\text{ is odd}.
\end{array}\right.\end{equation}
Thus equations (\ref{eqn4.2}) and (\ref{eqn4.5}) show $|\TKSp(M)|\le|\vec\eta_\nu(\mathcal{K}_\nu)|$. As the opposite
inequality is immediate, we have
$$\vec\eta_\nu(\mathcal{K}_\nu)=\mathcal{K}_\nu=\TKSp(M).$$
The Theorem now follows from equation (\ref{eqn4.1}).
\qedbox

\section{The groups $\Tko_{4k-1}(B\QL)$}\label{Sect-ko}
Let $x=(M,g,s,f)$ where $s$ is a spin structure and $f$ is a $G$ structure on a compact Riemannian manifold $(M,g)$ of dimension
$4k-1$.  Let $\eta^\sigma(x)$ be the eta invariant of the associated Dirac operator with coefficients in $f^*\sigma$. We reverse the
parities of the invariant defined in the previous section to define:
$$
\vec\eta_k(x):=\left\{\begin{array}{ll}
(\eta^{2\Theta_1}(x),\eta^{2\Theta_2}(x),\eta^{\Delta}(x),\eta^{2\Delta^2}(x)...,\eta^{2\Delta^k}(x))&(k\text{ even)}\\
(\eta^{\phantom{2}\Theta_1}(x),\eta^{\phantom{2}\Theta_2}(x),\eta^{2\Delta}(x),\eta^{\Delta^2}(x)...,\eta^{2\Delta^k}(x))&
(k\text{ odd}).
\end{array}\right.$$
We have used real representations if $k$ is odd and quaternion representations if $k$ is even. Therefore, by Corollary
\ref{corollary3.7},
$\vec\eta_k$ extends to:
$$\vec\eta_k:\Tko_{4k-1}(BG)\rightarrow(\Q/2\Z)^{k+2}.$$

The group $\QL$ has 3 non-conjugate elements of order $4$: $\{\mathcal{I},\mathcal{J},\xi\mathcal{J}\}$
which generate the $3$ non-conjugate
subgroups $\{\langle\mathcal{I}\rangle,\langle\mathcal{J}\rangle,\langle\xi\mathcal{J}\rangle\}$ of order
$4$. The representation $\gamma_1$ restricts to a fixed point free representation of any subgroup of $\QL$. We define the following
spherical space forms:
$$\begin{array}{ll}
M_{Q}^{4k-1}:=M^{4k-1}(\QL,k\gamma_1),&
M_{\mathcal{I}}^{4k-1}:=M^{4k-1}(\langle{\mathcal{I}}\rangle,k\gamma_1)\\
M_{\mathcal{J}}^{4k-1}:=M^{4k-1}(\langle{\mathcal{J}}\rangle,k\gamma_1),&
M_{\xi\mathcal{J}}^{4k-1}:=M^{4k-1}(\langle{\xi\mathcal{J}}\rangle,k\gamma_1).
\end{array}$$
Give the lens spaces $M_g^{4k-1}$ the $\QL$ structure induced by the natural inclusion $\langle g\rangle\subset\QL$. We project
into the reduced group
$\TMSpin_{4k-1}(\QL)$; this does not affect the eta invariant as $\eta^\sigma(\MSpin_*(pt))=0$. Let $i>0$.
By Theorem \ref{theorem3.8}:
$$
\begin{array}{rll}
(\eta^{\Theta_1},\eta^{\Theta_2},\eta^{\Delta^i})(M_{\mathcal{I}}^{4k-1}-M_{\mathcal{J}}^{4k-1})&=\left\{\begin{array}{llll}
2^{-k}(2,&1,&0)&\text{if }\ell=8,\\
2^{-k}(1,&0,&0)&\text{if }\ell>8,\end{array}\right.\\ \\
(\eta^{\Theta_1},\eta^{\Theta_2},\eta^{\Delta^i})(M_{\mathcal{I}}^{4k-1}-M_{\xi\mathcal{J}}^{4k-1})&=\left\{\begin{array}{llll}
2^{-k}(1,&2,&0)&\text{if }\ell=8,\\
2^{-k}(0,&1,&0)&\text{if }\ell>8,\end{array}\right.\\\\
(\eta^{\Theta_1},\eta^{\Theta_2},\eta^{\Delta^i})(M_Q^{4k-1})&=\begin{array}{llll}(0,&0,&c_{i-k})\phantom{....}&\text{any }\ell.
\end{array}
\end{array}$$
Let $K^4$ be a spin manifold with
$\hat A(K^4)=2$ and let
$B^8$ be a spin manifold with $\hat A(B^8)=1$. Let
$Z^{8k-4}:=K^4\times B^{8k-8}$ and $Z^{8k}=(B^{8})^k$. Standard product formulas \cite{PG89} then show
$$\eta^\sigma(M^{4k-1}\times Z^{4j})=\eta^\sigma(M^{4k-1})\hat A(Z^{4j})=\left\{
\begin{array}{rl}2\eta^\sigma(M^{4k-1})&\text{if }j\text{ is odd},\\
\eta^\sigma(M^{4k-1})&\text{if }j\text{ is even}.
\end{array}\right.$$
Let $B_\nu$ and $\mathcal{B}_\nu$ be as defined in Section \ref{Sect-KSp}. There is a dimension shift involved as we must set $\nu=k+1$. We
use the same arguments as those given previously to see
$$
\vec\eta_k\left(\begin{array}{l}
M_{\mathcal{I}}^{4k-1}-M_{\mathcal{J}}^{4k-1}\\
M_{\mathcal{I}}^{4k-1}-M_{\xi\mathcal{J}}^{4k-1}\\
M_Q^{4k-1}\\
M_Q^{4k-5}\times Z^4\\
...\\
M_Q^3\times Z^{4k-4}
\end{array}\right)=\left(\begin{array}{ll}C_k&0\\0&B_{k+1}
\end{array}\right)\in M_{k+2}(\Q/2\Z)$$
where $C_k$ is the $2\times 2$ matrix given by
$$
C_k=\left\{\begin{array}{ll}
2^{1-k\phantom{-1}}\left(\begin{array}{ll}
2&1\\1&2
\end{array}\right)&\text{ if }\ell=8\text{ and }k\text{ is even},\\
2^{1-k\phantom{-1}}\left(\begin{array}{ll}
1&0\\0&1
\end{array}\right)&\text{ if }\ell>8\text{ and }k\text{ is even},\\
2^{-k}\left(\begin{array}{ll}
2&1\\1&2
\end{array}\right)&\text{ if }\ell=8\text{ and }k\text{ is odd},\\
2^{-k}\left(\begin{array}{ll}
1&0\\0&1
\end{array}\right)&\text{ if }\ell>8\text{ and }k\text{ is odd},
\end{array}\right.$$

Theorem \ref{theorem1.2} will follow from Theorem \ref{theorem4.3} and from the following:

\begin{theorem}\label{theorem5.1} We have\rm
$$\Tko_{4k-1}(B\QL)=\left\{\begin{array}{ll}
\Z_{2^{k}}\oplus\Z_{2^{k}}\oplus\mathcal{B}_{k+1}&\text{if }k\text{ is even},\\
\Z_{2^{k+1}}\oplus\Z_{2^{k+1}}\oplus\mathcal{B}_{k+1}&\text{if }k\text{ is odd}.
\end{array}\right.$$\end{theorem}

\medbreak\noindent{\bf Proof:} We use  the same argument used to prove Theorem \ref{theorem4.3}. Let
\begin{eqnarray*}
&&\mathcal{L}_k:=\operatorname{Span}_\Z\{M_{\mathcal{I}}^{4k-1}-M_{\mathcal{J}}^{4k-1},
M_{\mathcal{I}}^{4k-1}-M_{\xi\mathcal{J}}^{4k-1},M_Q^{4k-1},\\
&&\phantom{\mathcal{L}_k:=\operatorname{Span_{\Z}}..}M_Q^{4k-5}\times Z^4,...,M_Q^3\times
Z^{4k-4}\}\subset\Tko_{4k-1}(B\QL).\end{eqnarray*}
 We then have that
$$\vec\eta_k(\mathcal{L}_k)=\left\{\begin{array}{ll}
\Z_{2^{k}}\oplus\Z_{2^{k}}\oplus\mathcal{B}_{k+1}&\text{if }k\text{ is even},\\
\Z_{2^{k+1}}\oplus\Z_{2^{k+1}}\oplus\mathcal{B}_{k+1}&\text{if }k\text{ is odd}.
\end{array}\right.$$
By Lemma \ref{lemma2.1} we have $c_0=\frac{\ell-1}\ell$ and thus $2c_0$ is an element of order $\ell$ in $\Q/2\Z$. We use the
diagonal nature of the matrix $B_{k+1}$ to see that:
$$|\vec\eta_k(\mathcal{L}_k)|\ge\left\{\begin{array}{ll}
4^{k}\ell^k&\text{if }k\text{ is even},\\
4^{k+1}\ell^k&\text{if }k\text{ is odd}.\end{array}\right.$$ 
We use \cite{BGS97} Theorem 2.4 see:
$$|\Tko_{4k-1}(B\QL)|=\left\{\begin{array}{ll}
4^{k}\ell^k&\text{if }k\text{ is even}.\\
4^{k+1}\ell^k&\text{if }k\text{ is odd.\qquad\qedbox}\end{array}\right.$$

\begin{remark}\label{remark5.2}\rm Let $n\ge0$. One has \cite{BaBr96} that:
$$\begin{array}{ll}
\Tko_{8n+\varepsilon}(\Sigma^{-1}BS^3/BN)=&\left\{\begin{array}{ll}
\Z_2\phantom{...............a}&\text{if }\varepsilon=1,2,\\
\Z_{2^{2n+2}}&\text{if }\varepsilon=3,7,\\
0&\text{if }\varepsilon=4,5,6,8,\end{array}\right.\end{array}$$
We may use equation (\ref{eqn4.6}) to decompose:
\begin{eqnarray*}\Tko_*(B\QL)&=&
\Tko_*(\Sigma^{-1}BS^3/BN)\oplus\Tko_*(\Sigma^{-1}BS^3/BN)\\&
\oplus&\Tko_*(BSL_2(\mathbb{F}_q)).\end{eqnarray*}
This is the decomposition given in Theorems \ref{theorem4.3} and \ref{theorem5.1}:
\begin{eqnarray*}
\mathcal{A}_k&=&\Tko_{4k-1}(\Sigma^{-1}BS^3/BN)\oplus\Tko_{4k-1}(\Sigma^{-1}BS^3/BN)\\
&=&\operatorname{Span}\{[V^{\Theta_1}],[V^{\Theta_2}]\}\subset\TKSp(M^{4k+3}(\QL,\tau))\\
&=&\operatorname{Span}\{[M^{4k-1}_{\mathcal{I}}-M^{4k-1}_{\mathcal{J}}],
[M^{4k-1}_{\mathcal{I}}-M^{4k-1}_{\xi\mathcal{J}}]\}\subset\Tko_{4k-1}(B\QL),\\
\mathcal{B}_k&=&\Tko_{4k-1}(BSL_2(\mathbb{F}_q))\\
&=&\operatorname{Span}_\Z\{[V^{\varepsilon_j\Delta^j}]\}\subset\TKSp(M^{4k+3}(\QL,\tau))\\
&=&\operatorname{Span}_\Z\{[M_Q^{4k-1-4\mu}\times Z^{4\mu}]\}\subset\Tko_{4k-1}(B\QL).
\end{eqnarray*}
\end{remark}

\section*{Acknowledgements} Research of P. Gilkey partially supported by the NSF (USA) and the MPI (Leipzig, Germany).

\end{document}